\documentclass[12pt, twoside,a4paper]{article}
\pdfoutput=1
\usepackage[cp1251]{inputenc}
\usepackage{ifthen}
\usepackage{amsmath}
\usepackage{amsfonts}
\usepackage{latexsym}
\usepackage{amsthm}
\usepackage{amssymb}
\newcommand{\CopyName}{ V.\ M.\ Zhuravlov}
\newcommand{\NAME}{ V.\ M.\ Zhuravlov}
\newcommand{\Year}{2023}
\newcommand{\rightheadtext}{Uniformity and nonuniformity}
     \newcounter{chapter}
     \newcounter{artpage}[chapter]
     \newcommand{\vs}{\vspace{.1in}}
     \newcommand{\vsk}{\vspace{.2in}}
\makeatletter
     \renewcommand{\@evenhead}{\footnotesize \ifthenelse{\value{artpage}=0}
     {\hfil}{\thepage\hfil \textsc {\leftmark} \hfil } }
     \renewcommand{\@oddhead}{\footnotesize\ifthenelse{\value{artpage}=0}
     {\hfil}{\hfil \textsc \rightmark \hfil \thepage} }
     \newcommand{\logo}{\baselineskip2pc \hbox to\hsize{\hfil\copyright\,\footnotesize
     \CopyName, \Year}}
     \renewcommand{\@oddfoot}{\ifthenelse{\value{artpage}=0}{\logo
     \refstepcounter{artpage}} {\hfil\refstepcounter{artpage}}}
     \renewcommand{\@evenfoot}{\ifthenelse{\value{artpage}=0}{\logo
     \refstepcounter{artpage}} {\hfil\refstepcounter{artpage}}}
     \renewcommand{\section}{\@startsection{section}{1}{0pt}{3.5ex plus
     1ex minus .2ex}{2.3ex plus 2.ex}{\large\hfil\textsc}}
     \vspace{3pt}\vs

\vs
\newcommand{\tit}{Uniformity and nonuniformity}
\date{2023}
\usepackage{hyperref}
\begin{document}
\hfill
\vspace{0.3in}
\markboth{{\NAME}}{{\rightheadtext}}\begin{center} \textsc {\CopyName} \end{center}\begin{center} \renewcommand{\baselinestretch}{1.3}\bf {\tit} \end{center}
\vspace{20pt plus 0.5pt} {\abstract{\noindent
The article introduces the concept of uniformity, which is formulated as a scheme of axioms. The connection of this concept with ordered sets is studied. The effectiveness of using axiom schemes as a convenient and short way of replacing axiom systems is shown. With such a replacement, the main idea of the system of axioms of the given mathematical theory immediately becomes clear. The reasons for using number systems in mathematics are also shown.\newline
\textit{Uniformity and nonuniformity, 2023, udc: 510.65 msc: 03B16\vspace{3pt}}\newline
}
}\vsk
\tableofcontents
\section{Introduction: Relationship. Similarity and distinction}
We will look for the meaning of known mathematical concepts and try to create new concepts that can reveal other meanings...
In all logical formulas, we believe that unary connectives (for example, negation) are stronger than binary ones, conjunction and disjunction are equivalent, but stronger than implication, and logical equivalence is the weakest. In non-obvious cases, we use brackets; "$\equiv$" means identity or definition; curly braces mean a set of elements described inside the brackets.
Unless otherwise stated, for the time being we assume that we are in the universe of traditional set theory. At the same time, we will study the properties of the uniformeity of the structure of objects on which certain predicates are defined, without yet going into the specifics of these predicates, for example, we do not consider functions and functional terms, we will deal with this interesting topic later. We will try to ensure that everything that we express in words has the maximum possibility of formalization (and vice versa).\par
\textit{Any statement can be reduced to the form: "some objects distinguishable by means of \textbf{R} are identical by means of \textbf{E}".}\par
In other words, something different in one sense is similar in another sense. It seems that each statement can be represented as a pair \textbf{(\textit{R, E})}, where \textbf{\textit{R}} is some (quite detailed, detailed) relation of difference, and \textbf{\textit{E}} is some equivalence. We understand distinguishability as the negation of some equivalence (belonging to its various classes) and give its positive definition. Thus, one can speak of a pair of equivalences (together with the negation operation).
It is known that if, from a classical point of view, we look at the equivalence of \textbf{\textit{E}} on a certain set \textbf{\textit{M}}, then it will be: $E\subseteq(M\times M)$, $(\triangle_{M}\subseteq E)$, $(E^{-1}=E)$, $(E^{2}\subseteq E)$,\newline
— we define: $E^{2}\equiv {(x,z)\mid\exists y: (x,y)\wedge (x,z)}$ — we deliberately wrote out all these formulas in detail in order to determine distinguishability, as an addition to equivalence:
$$R(E)=(M\times M)-E(R)$$
— we have chosen a functional dependence to emphasize the fact that any way of distinguishing is always associated with a certain way of identification, and vice versa. Hence:
$$R\cap \triangle_{M}=\emptyset$$
$$R=R^{-1}$$
$$R\subseteq M^{2}-(M^{2}-R)^{2}$$
The first two equalities mean obvious irreflexivity: $\neg R(x,x)$, and symmetry: $R(x,y)\Leftrightarrow R(y,x)$. As for inclusion, according to the transitivity of equivalence:
$$((x,y)\notin R)\wedge((y,z)\notin R)\Rightarrow((x,z)\notin R)$$
— or, in other words:
$$\forall(x,y,z):xRz\Rightarrow((yRx)\vee(yRz))$$
It can also be output like this:
$$(E(x,y)\wedge E(y,z))\Rightarrow E(x,z)$$
$$\neg E(x,y)\vee\neg E(y,z)\vee\neg E(x,z)$$
$$\neg E(x,z)\Rightarrow(\neg E(x,y)\vee \neg E(y,z))$$
where \textbf{\textit{E}} is the equivalence corresponding to \textbf{\textit{R}}. Thus, the distinguishability of \textbf{\textit{R}} is an irreflexive, symmetric, and antitransitive relation. Anti-transitivity (negation of transitivity) means that every element differs from one of two distinct elements...
Here is a list of the abstract binary relations we use:\newline\par
\begin{tabular} { |c|c|c|c| }
\hline
equivalence&1&1&1\\
\hline
distinguishability&1&0&-1\\
\hline
orderliness&-1&1&1\\
\hline
strict order&-1&0&1\\
\hline
preorder& &1&1\\
\hline
 &symmetry&reflexivity&transitivity\\
\hline
\end{tabular}\newline\newline
— where: \textbf{1} is true, \textbf{0} is negative true, \textbf{-1} is antisymmetric, antitransitive. Not very necessary 4-th property is to add linearity (2-completeness): $(a\neq b)\Rightarrow ((a>b)\vee(b>a))$. We will also use one ternary relation — the cyclic order:\newline\par
\begin{tabular} { |c|c|c|c| }
\hline
cyclic order&$[a,b,c]\Rightarrow \neg[c,b,a]$&$[a,b,c]\wedge[a,c,d]\Rightarrow[a,b,d]$&$[a,b,c]\Rightarrow[b,c,a]$\\
&3-asymmetry&3-transitivity&cyclicity\\
\hline
\end{tabular}\newline\newline
— and linearity (3-completeness) is also possible for it: $[b\neq c]\wedge[c\neq a]\Rightarrow([a,b,c]\vee[c,b,a])$;
it should also be clarified here that the non-execution of a property does not mean the fulfillment of its negation: the truth or falsity of a property must be satisfied on all pairs from \textbf{\textit{M}}; therefore, one should distinguish between irreflexivity and non-reflexivity, anti-transitivity and non-transitivity, etc…\par
For order relations, the properties of density and continuum are also important.\par
Interestingly, all complete relations (when all sequences of elements are comparable "from the inside") have limitingly complete closures to equivalence $(M\times M)$ and to distinguishability $(M\times M)-\triangle_{M}$, i.e. — give a certain instrumental opportunity for maximum assimilation or differentiation of elements.
All the relations in the list are so widely used in mathematics because they help the effective use of similarity and difference. For the time being, we remain within the framework of set theory, and therefore we consider binary and ternary relations to be basic.\par
It is clear that ordering implies trivial distinguishability (as the negation of minimal equivalence-equality). We believe that the basic properties and relations must always be determined effectively, and this efficiency, in turn, must itself be relativized with respect to various kinds of "means". This means that distinguishability, equivalence, similarity, equality, etc. — must have their basis in something else... A good means of distinguishing and identifying is, for example, the relation of order, as well as systems of certain algebraic operations.\par
A few words about the goals of our constructions. The axiom of choice proclaims the possibility of well ordering (in the sense that each subset of a set of individuals has a lower bound) of any set, which causes great objections among constructivists. We believe that constructive constructions are strongly dependent on means, which are not always the same, absolute and immutable; formalism rests on the other extreme, believing that mathematics abstractly owns all means in general. It seems that in any set, one should first understand what we distinguish in it, where equivalence becomes what is called trivial equality? Having distinguishable elements, we can create a transitive closure of any collection of pairs (preliminarily cutting it down to an irreflexive and antisymmetric collection),— then we get a linear order (do not forget that "total" and "well" are different things; total indicates comparability of any elements, the well order indicates the presence of the lower bound of any set). Of course, in the finite case, this procedure is easily algorithmized, but it requires recursion (and even transfinite) in the infinite case. So perhaps we should accept the following, weakened axiom of choice:\newline\par
\textit{Any distinguishable set can be \textbf{linearly} ordered.}\newline\par
Here we have so far only paid attention to the importance of likening and distinguishing objects and have given a positive formalization of the relation of distinguishability (which was previously intuitively understood as the negation of equivalence). Further study of this issue will require slightly different approaches.
\section{UNIFORM THEORIES}
As Bertrand Russell stated: "...mathematical knowledge has ceased to be mysterious. It has the same nature as the "great truth", which is 3 feet in a yard." Later, however, it turned out that this is not entirely true with mathematics, and with numbers it is also not so at all. In addition to numbers and trivial tautologies, mathematics studies things that are much more abstract, non-obvious and even incalculable. And yet, a huge amount of truth in the words of Russell is undeniable to this day. Indeed, why is arithmetic the necessary foundation for all mathematical (and indeed linguistic) constructions? And when we lay the foundation on other foundations (sets, categories, algebras…), we still cannot bypass arithmetic… Why is almost all exact science based on measurement — an empirical reflection of nature in a certain numerical system, one way or another generalizing all the same arithmetic? After all, there are numbers all around, from all sorts of alphabets and formal theories, to cars and astronomical observations. Let's try to find an answer. All our cognitive activity can be interpreted on the basis of the operations of assimilation and distinction. And these operations are the more productive, the more "powerful" a certain universe is, potentially capable of accommodating something maximally similar and different at the same time. Thus we arrive at a set of uniformities objects, under certain conditions, capable of becoming as heterogeneous and distinct as possible. Let's start with the formalization of the concept of uniformity.\par
\textit{We already know that classical distinguishability is a binary, irreflexive, symmetric, and antitransitive relation.} Let's try to define logical distinctness.\par
Any unary formula (it can also contain constants, and its other variables can be connected by quantifiers) can serve as an indicator for identifying and distinguishing the individual elements on which it is defined. The corresponding equivalence and distinguishability are binary relations: $A(x)\Leftrightarrow A(y)$ and its negation $\neg[A(x)\Leftrightarrow (y)]$; thus, each solvable formula will define such distinguishability (and equivalence) on each model set; and every classical theory induces its own trivial equality by means of all its indicators. And outside of its models, the theory distinguishes those models from the objects on which it is not true. \textit{Therefore, any constant-free unary formula will be called an indicator.} The requirement of being constant follows from the fact that constants are distinguishable a priori. Considering constants, we can also talk about an indication, but this will already be an indication localized with respect to a given constant. We will return to this discussion later.\par
The set of all distinguishability on the model set M is an inclusion lattice (it is dual to an equivalence lattice). The set of indicators of the theory forms a lattice with respect to the implication. This lattice is isomorphic to some sublattice of distinguishability on \textbf{\textit{M}}. Such a sublattice also includes empty distinguishability, since any theory also has zero indicators. \textbf{1}\textit{-uniformity is equivalent to the fact that all indicators of the theory are zero.} This means that the distinguishability of the elements of \textbf{\textit{M}} can only be achieved a priori by introducing constants. In other words, all elements of \textbf{\textit{M}} are symmetrical. But if \textbf{\textit{M}} is one-element, then it is both finitely uniformity and also models any finitely uniformity theory. However, its information capacity is minimal. Therefore, it should be acknowledged substantively that of the two models of a uniformity theory, the model of higher power will be more uniformity.\newline\par
\textbf{Definition 1} \textit{We call a classical formal theory \textbf{1}-uniformity if:
$$\exists x:A(x)\Rightarrow \forall x:A(x)$$
for any non-constant unary formula \textbf{A(x)} with a single free variable; in other words, it is the scheme of the axioms of the theory. A linguistically invariant extension (see the book by Rasiowa and Sikorsky) \textbf{T2} of the theory \textbf{T1} will be called more uniformity than \textbf{T} if any formulas of the type:
$$\bot...\exists x:A(x)\bot...\Longrightarrow \bot...\forall x:A(x)\bot...$$
(where: $\bot$ — any quantifier,— and all such quantifiers with their variables coincide on the right and left of the implication; \textbf{A} — any formula) are fulfilled in \textbf{T2} if they are fulfilled in \textbf{T1}; sets that model such schemes of uniformity in an obvious way are also called uniformities.} Thus, the set of linguistically identical theories (not strictly! — in our definition) is ordered by uniformity.\newline\par
We see that a strongly simple (that is, trivial) theory is \textbf{1}-uniformity, and in general, maximally uniformity, since all variables of any formula can be related by the quantifier $\forall$. Thus, the degree of uniformity characterizes gradations of strong simplicity. The more uniformity the theory, the more difficult it is to distinguish its objects, because the field of indicators narrows. The requirement of non-constancy follows from the fact that the introduction of any constant means a priori \textit{non-uniformity}, distinguishability. If we switch to open formulas in the definition, then it will be: $A(c)\Rightarrow A(x)$ where c is a constant.\par
Constants can be introduced into any theory. In uniformity — too. And then for the constant c we have: $\exists x:(x=c)\wedge\exists x:(x\neq c)$, which already entails non-uniformity. In this case, all individuals will still be a uniformity set, because we are still dealing with the same theory. That is why in the definition of uniformity we require the non-constancy of indicators. And then the fixation of a certain constant sets the "induced" non-uniformity, localized with respect to this constant. This is the value of a uniformity set — informationally empty, but capacious, capable of containing a lot of information.\par
However, non-uniformity can arise not only from constants. The axiom: $\exists x\forall y:(x\ast y =y)$ of a constant group theory is essentially non-uniformity, unless the group is nontrivial. Any theory can be written only in terms of predicate formulaic symbols, without functors and constants. Then the existential statement about the presence of a certain constant will be: $\exists!x...$. A theory whose constants are determined only in this way will be called existentially non-uniformity (more precisely, \textit{1-non-uniformity}). Moreover, it can also be \textit{1-uniformity} (as the theory of unbounded linear order — see below).\newline\par
\textbf{Definition 2} \textit{We call \textbf{\textit{T} n}-uniformity if the theory induced by it on variables spanning all \textbf{n}-element subsets of the set spanned by \textbf{T} is \textbf{1}-uniformity. A theory that is \textbf{n}-uniformity for all finite \textbf{n} is said to be finitely uniformity.}\newline\par
As for \textbf{1}-uniformity: $\exists x:A(x)\Rightarrow \forall x:A(x)$, it will be interesting to look at its scheme in a different, equivalent form:
$$\forall x:A(x)\bigvee \forall x:\neg A(x)$$
— according to the above, we mean by \textit{\textbf{A(x)}} any, at most, unary-open formula. Therefore, for example, any binary relation on a 1-uniformity set must be either everywhere reflexive or else everywhere irreflexive; similarly, \textbf{2}-uniformity requires a binary relation to be symmetric or antisymmetric — on all pairs in the set; in \textbf{3}-uniformities sets, the existence of at least one transitive triple of distinct elements requires the transitivity of the entire binary relation (moreover, the existence of at least one pair transitive for any third element requires the same transitivity even in the case of \textbf{2}-uniformity).\par
Formulas, even more open, are equivalent to the same ones, but given on some Cartesian powers of an individual domain. The latter are not very uniformities simply by virtue of the definition of the Cartesian degree. Namely, for any set \textit{\textbf{M}}, $diag(M^{2})$ and $M^{2}-diag(M^{2}$ are distinguishable; the set $diag(M^{2})$ has symmetry in \textit{\textbf{(x,y)}} and \textit{\textbf{(y,x)}}. The set $diag(M^{2})$ is isomorphic to \textit{\textbf{M}}, that is,— its internal properties are in the "department" of \textbf{1}-uniformity... In other words, to talk about the properties of two-element subsets, one must consider $M^{2}-diag(M^{2})$ without structure. So, in the scheme of \textbf{2}-uniformity, we need to deduce from the existence of exactly a two-element unordered pair its universality. \textit{Therefore, the \textbf{2}-uniformity scheme will be:}\par
$$\exists x\exists y:(x\neq y)\wedge A(x,y)\Rightarrow \forall x\forall y:(x=y)\vee A(x,y)\vee A(y,x)$$
\textit{Or, (as in the case of \textbf{1}-uniformity) in a different form:}
$$\forall x\forall y:(x=y)\vee A(x,y)\vee A(y,x)\bigvee\forall x\forall y:\neg A(x,y)\vee\neg A(y,x)$$
\textit{where \textbf{A(x, y)} is a non-constant binary formula; all such variants of uniformity define these schemes in a disjunctive (and in practice, in a preliminary!) form.
Similarly, we define \textbf{3}-uniformity:}
$$\exists x\exists y\exists z:(x\neq y)\wedge(y\neq z)\wedge(z\neq x)\wedge A(x,y,z)\Longrightarrow\forall x\forall y\forall z:(x=y)\vee(y=z)\vee(z=y)\vee$$
$$\vee A(x,y,z)\vee A(y,x,z)\vee A(x,z,y)\vee A(y,z,x)\vee A(z,x,y)\vee A(z,y,x)$$
\textit{which in disjunctive form:}
$$\forall x\forall y\forall z:(x=y)\vee(y=z)\vee(z=x)\vee\neg A(x,y,z)\bigvee$$
$$\bigvee\forall x\forall y\forall z:(x=y)\vee(y=z)\vee(z=x)\vee A!(x,y,z)$$
\textit{where the "factorial truth" of $A!(x_{1},…,x_{n})$ means that there is a permutation in the sequence $(x_{1},…,x_{n})$ for which \textbf{A} is true.}\par
\textit{In general, the scheme of the \textbf{n}-uniformity axioms will be:}
$$\exists(x_{1},…,x_{n}):\bigwedge_{i,j\in(1...n)}(x_{i}\neq x_{j})\wedge A(x_{1},…,x_{n})\Longrightarrow\forall(x_{1},…,x_{n}):A!(x_{1},…,x_{n})$$
\textit{where the quantifier $\bot(x_{1},…,x_{n})$ means $\bot x_{1}...\bot x_{n}$; disjunctively:}
$$\forall(x_{1},…,x_{n}):\bigvee_{i,j\in(1...n)}(x_{i}=x_{j})\vee\neg A(x_{1},…,x_{n})\bigvee$$
$$\bigvee\forall(x_{1},…,x_{n}):\bigvee_{i,j\in(1...n)}(x_{i}=x_{j})\vee A!(x_{1},…,x_{n})$$
Verbally, this all means: if a non-constant \textit{\textbf{n}}-open formula \textit{\textbf{A}} is true for some \textit{\textbf{n}} different variables, then it is true factorially for all such \textit{\textbf{n}}. Or, in other words, if \textit{\textbf{A!}} is true for one \textit{\textbf{n}} different variables, then it is true for all such \textit{\textbf{n}}. Thus, we have formalized the previously given meaningful concept of \textit{\textbf{n}}-uniformity for any finite \textit{\textbf{n}}.\par
These are all axiom schemas, so they are no less powerful in application than schemas from set theory, such as the principle of induction or the substitution schema. The concept of uniformity is applicable to the theory \textit{\textbf{T}}, to the set \textit{\textbf{M}} that models such a theory, as well as to a specific formula \textit{\textbf{A}}, which can be: a statement scheme (as in our definitions) or some kind of predicate. Those — we can talk about theories, individuals, or about predicates. From the context, unless otherwise stated, it will always be clear to what we refer this concept. It is clear that:\par
\textit{If certain formulas are more or less uniformies, then so will be their conjunctions, disjunctions, and negations. But I don't know how uniformity behaves when quantifying some variables in such formulas?}\par
All these schemes are obtained based on the fact that the properties of any diagonals of the set $M^{n}$ are reduced to the properties of $M^{n-k}$, and also on the fact that uniformity must be invariant under transformations that preserve the Cartesian structure of $M^{n}$. In short, the \textit{\textbf{n}}-uniformity of \textit{\textbf{T}} is \textbf{1}-uniformity for the theory induced by our \textit{\textbf{T}} on the set of all \textit{\textbf{n}}-element subsets of an individual domain, according to the above meaningful definition. We have come to this construction because the sequences of elements are our way of looking at them, and not the properties of \textit{\textbf{T}} itself, and only by "twisting" these sequences as much as possible will we obtain the properties of \textit{\textbf{T}}.
\section{Powers of models and propositions}
We are talking about the power of a statement in the sense that any model of it can have its own power of its truth area.\newline\par
\textbf{Definition 3} \textit{A theory \textbf{T} is called (n,m,l,…)-predicate if the numbers \textbf{n,m,l,…} exhaust the arity of all predicates that are non-autologically included in the theorems of \textbf{T}, expressed exclusively predicately (without functors and constants). We note that expressing \textbf{T} in terms of predicates of other arities, we get another \textbf{T}. But their structures can, in a certain way, be attributed to each other.}\newline\par
Now some facts about n-uniformities \textit{\textbf{(n,m,l,…)}}-predicate \textit{\textbf{T}}. Any \textbf{1}-uniformity theory does not distinguish its constants introduced from outside, a priori. These are theories about what would happen if there were a means of absolute discrimination. All elements of the set described by a \textbf{1}-uniformity theory are symmetric, if we do not consider its constants. This is, as it were, the simplest, linear form of distinguishability. A strongly simple theory is generally unable to distinguish between its constants, for it they are all identical. The theory of trivial equality is also \textit{\textbf{n}}-uniformity for any, even infinite, \textit{\textbf{n}}. Further, according to the \textbf{1}-uniformity scheme, it will be:
$$\exists x\forall y:P(x,y)\Rightarrow\forall x\forall y:P(x,y)$$
— which leads either to triviality or to a contradiction; hence, such \textit{\textbf{x}} does not exist in non-trivial consistent \textbf{1}-uniformities theories of any binary predicate \textit{\textbf{P}}, i.e. — in the case of an order relation, the individual area is not limited from above and below and is infinite.\par
Looking at the key properties of binary relations from the above table, we see that reflexivity must be considered a \textbf{1}-uniformity predicate (and in general, of course, uniformity, because it is unary, and therefore we will supplement it with any sequence of dummy variables); so will the symmetry (or antisymmetry, whichever it is) with respect to \textbf{2}-uniformity (or higher); the same is true for (binary) transitivity (\textbf{3}-and more-uniformity). And this will be true for any Boolean constructions from these statements, as we have already noted. Thus we have proved:\newline\par
\textbf{Lemma 1} \textit{A non-trivial, linearly ordered set is \textbf{1}-uniformity if and only if it has no least and no largest elements. This set will be finitely uniformity (\textbf{n}-uniformity for all finite \textbf{n})\nopagebreak[4] if and only if it is everywhere dense.}\newline
It is in this case that pairs of elements are indistinguishable by the number of elements lying between them.\newline\par
As for the other relations we have defined, then a non-trivial complete equivalence whose classes are equivalent in power will always be \textbf{1}-uniformity, and will never be \textbf{2}-uniformity, due to the fact that there are equivalent and non-equivalent pairs. The same is true of distinguishability, for it is the negation of equivalence. And the same should be said about such "partial" relations as preorder and comparability. Similarly, direct products and direct sums of uniformities sets are always \textbf{1}-uniformity, but not \textbf{2}-uniformity. We have already written about the order, but it is worth adding that well-ordered sets are "maximally" non-uniformities; moreover, among those, the sets of topics are more non-uniformities, the more powerful they are, since a larger number of different elements makes it possible to "diversify" their indicators more. \textit{Therefore, the natural series is more non-uniformity than finite sets, and the ordinal series is even more non-uniformity. Further, the ternary relation of a cyclic order is, in essence, an order localized with respect to the elements of a set. And if it is linear, then it is \textbf{1}-uniformity, and if it is everywhere dense, then it is \textbf{n}-uniformity,} just like the usual ordering. It is interesting to note that \textit{only the cyclic order makes it possible to "rebuild" a completely non-uniformity finite set (that is, a well-ordered one) into a \textbf{1}-uniformity one.}\par
Let's talk about the power of individual regions. Theories with an empty individual domain are controversial. If the region \textit{\textbf{M}} consists of one element,— $\mid M\mid=1$, then its theory \textit{\textbf{T}} will be absolutely homogeneous for any \textit{\textbf{n}}. When $\mid M\mid=2$, its \textit{\textbf{T}} will be \textbf{2}-uniformity. Let $\mid M\mid>2$, and its theory is \textbf{2}-uniformity. Let \textit{\textbf{R(x)}} be any constant-free \textbf{1}-open formula. Substituting the formula $B(x)\wedge \neg B(y)$ into the \textbf{2}-uniformity scheme, we obtain:
$$\forall x:B(x)\vee\forall x:\neg B(x)$$
Consequently, a \textbf{2}-uniformity \textit{\textbf{T}} with an individual domain $\mid M\mid>2$ (more precisely, we mean not even the power of the model, but the truth of formulas that assert the existence of more than two different elements) will also be \textbf{1}-uniformity. Hence we conclude that if the theory is \textbf{2},— but not \textbf{1}-uniformity, then $\mid M\mid=2$. Similarly, it is easy to prove (if there are only \textbf{\textit{n}} elements, then there is no sequence of \textbf{\textit{(n + 1)}} pairwise unequal elements):\newline\par
\textit{\textbf{Lemma 2}: If for a theory \textbf{T} it is true that $\mid M\mid=n$, then for all \textbf{m}, \textbf{T} will also be \textbf{(n+m)}-uniformity.}\par
The lemma can be rephrased:\par
\textit{If \textbf{T} is n-uniformity for $n>1$, then $\mid M\mid>n$.}\par
And further:\par
\textit{The scheme of \textbf{n}-uniformity is the $(1...(n-1))$-identity formula.}\par
In short, the \textit{\textbf{n}-uniformity condition is a $(1...(n-1))$-tautology.} This can be verified by simulating this condition in an \textit{\textbf{n}}-element individual domain.\newline\par
We also add that if from the restrictions on the power of the model $\mid M_{T}\mid$ some restrictions on the nature of uniformity follow, this does not mean that the opposite is true. Thus, \textit{\textbf{n}}-identical formulas can be executed in an individual domain of a different, even infinite cardinality (and even more so, under restrictions on substitutions in circuits). The theory of unbounded linear order is \textbf{1}-uniformity, while the scheme of \textbf{1}-uniformity is \textbf{1}-identical, while $\mid M_{T}\mid$ at least countable.\par
It should also be noted that n-uniformity imposes some restrictions on statements of less than textit{\textbf{n}} dimensions. So, any unary non-constant \textit{\textbf{A(x)}} forms an textit{\textbf{n}}-dimensional conjunction: $\bigwedge_{i=1...n}A(x_{i})$ on sequences of different $x_{i}$. But then: If property \textit{\textbf{A}} is fulfilled on \textit{\textbf{n}} different elements of an \textit{\textbf{n}}-uniformity \textit{\textbf{M}}, then it is fulfilled on the whole of \textit{\textbf{M}}. Similarly, \textit{if \textbf{A} is fulfilled on one \textbf{M}-element, then it will be fulfilled on one of any \textbf{n} elements} (we find out this, substituting the disjunction into the scheme of \textit{\textbf{n}}-uniformity: $\bigvee_{i=1...n}A(x_{i})$).\par
Let our \textit{\textbf{T}} be \textit{\textbf{n}}-uniformity, but not \textbf{1}-uniformity. Hence, it contains a non-constant sentence: $B(x)\neq0,1$; let's assume that it is solvable. Then:
$$0<|\{x\mid B(x)=1\}|<1$$
$$0<|\{x\mid B(x)=0\}|<1$$\par
Means:
$$1<|M|<2n$$\par
As you can see, the distribution of uniformity on the natural axis significantly limits the power of the $M_{T}$.\newline\par
\textbf{Theorem 1}: \textit{If \textbf{T} is \textbf{(1...n)}-uniformity and \textbf{(1...n)}-predicate, then it is finitely uniformity.}\newline\par
Proof. Let we have a non-constant fully open formula:
$$A_{m}(y_{1}...y_{m})\wedge B_{k}(t_{1}...t_{k})$$
where: $A_{m}$ and $B_{k}$ are any formulas $(k,m < n+1)$ without constants; $y_{i}, t_{h}$ are their free variables; throughout the derivation it is assumed that: $y_{i}, y_{h}$, as well as $t_{i}, t_{h}$ are different for different \textit{\textbf{i, h}}; Then:
$$[\exists (y_{i\in 1...m}):A_{m}(y_{i\in 1...m})\wedge\exists (t_{j\in 1...k}):B_{k}(t_{j\in 1...k})]\Rightarrow$$
$$\Rightarrow[\forall (y_{i\in 1...m}):A_{m}!(y_{i\in 1...m})\wedge\forall (t_{j\in 1...k}):B_{k}!(t_{j\in 1...k})]\Rightarrow$$
$$\Rightarrow[\forall (y_{i\in 1...m},t_{j\in 1...k}):(A_{m}(y_{i\in 1...m})\wedge B_{k}(t_{j\in 1...k}))!]$$
We have shown the \textit{\textbf{(m+k)}}-uniformity of the original formula. The \textit{\textbf{(m+k+…)}}-uniformity of any chain of such $A_{m}, B_{k},C_{l}$, connected by conjunction or disjunction is shown in the same way. But in this way, all non-constant formulas \textit{\textbf{T}} have been exhausted. For if the original formula is not completely open, i.e., there is a quantifier prefix in front of it, knitting variables identical in $A_{m}$ and $B_{k}$, then we will rename, replace all quantifiers from the prefix with "$\exists$", we will carry out the previous proof, and as a consequence, we will prove the uniformity of the original formula. If in the formulas $A_{m}$ and $B_{k}$ some free variables coincide, then replacing them with non-coinciding ones, we carry out the same proof. And we pay attention to the fact that in the resulting formula \textit{\textbf{(m+k)}}-uniformity there are no conditions for the inequality of the variables $y_{i}$ and $y_{h}$; this means that if some $(t_{i} = y_{h})$, then even then there is such a permutation $(t_{h})$ in $B_{k}$ and a permutation $(y_{i})$ in $A_{m}$ that the formula $A_{m}\wedge B_{k}$ will be "as much"-uniformity as there are different free variables in it. Taking into account the previous discussion about quantifier prefixes, we obtain that any constant-free formula will be uniformity (i.e., we will exhaust all formulas). The theorem has been proven.\newline\par
It is also quite obvious that:\par
\textit{Finite uniformity unary predicate theories are strongly simple.}
\section{Uniformity and order}
\textbf{Theorem 2}: \textit{If a theory is finitely uniformity, binary, not very simple, and nontrivial, then it is a theory of everywhere dense unbounded linear order with a single predicate. The reverse is also true.}\newline\par
Proof. By nontriviality we mean nonidentity with the theory of trivial equality. So the axioms of \textit{\textbf{T}} must somehow bind the \textbf{2}-predicates. If \textit{\textbf{P}} is one of them, then \textbf{1}-uniformity gives:
$$\forall x: P(x,x)\vee\forall x:\neg P(xx)$$
— the obvious uniformity of the unary predicate;
$$\nexists x\forall y:P(x,y)\wedge\nexists x\forall y:P(y,x)$$
— this follows from the denial of strong simplicity; and further:
$$\forall x\exists y:P(x,y)\wedge\forall x\exists y:P(y,x)$$
— since non-triviality and not strong simplicity require that: $\exists x\exists y:P(x,y)$.
\textbf{2}-uniformity now gives:
$$\forall x\forall y:P(x,y)\vee P(y,x)$$
$$\forall x\forall y:(x\neq y)\Rightarrow\neg(P(x,y)\wedge P(y,x))$$
$$\forall x\forall y\exists t: P(x,y)\Rightarrow (P(x,t)\wedge P(t,y))$$
(the pairwise inequality of the entire triple is implied; the existence of a "condensed" pair follows from the previous one). The first formula again follows from the existence of a comparable pair; and from this, and also from the fact that all pairs are either symmetric or antisymmetric, we obtain exactly the antisymmetry (ie, the second formula), otherwise we would get the trivial theory of the total relation.\par
\textbf{3}-uniformity gives us transitivity, since for different \textit{\textbf{x, y, t}} the existence of a transitive triple follows from the previous formula (a transitive triple, for example, will be the one obtained by condensing a pair), from which we deduce that all triples are transitive.\par
We have proved that \textit{\textbf{P(x, y)}} is an (ordinary or strict) everywhere dense unbounded linear ordering. The theory of such an ordering is maximal, so no other properties can be added to it. We proved the converse theorem earlier. Note that the \textbf{2}-order predicate, once it is given, becomes unique, since the existence of a pair of elements for which one of the possible predicates is true and the other is false implies that the same is true for all pairs, i.e. — any other ordering will be opposite to the first one (i.e. — due to uniformity, we can choose an arbitrary ordering of the model set, but the choice is always alternative: there will be only one ordering). The theorem has been proven.\par
We see that the property of everywhere density plays an important role for finite uniformity; the same order, but discrete, will only be \textbf{1}-uniformity, because there will be pairs that cannot be compacted (and also pairs that can be compacted a finite number of times). It is interesting that for an everywhere dense order it will be: $P^{2}=P$, while for a discrete one we have a strict inclusion: $P^{2}\subset P$. The axiom of transitivity is implicative, and in the discrete case for uncondensed pairs it is purely formal; while for density everywhere, each of the pairs will be transitive in the sense of material implication (the element lying between the pair will form a transitive triple with it).
\section{Uniformity and cyclical order}
Let us now consider the ternary relation of cyclic order in more detail. Intuitively, such a relation should be understood as: "\textit{\textbf{b}} lies between \textit{\textbf{a}} and \textit{\textbf{c}}" or "after a, first comes \textit{\textbf{b}}, and then \textit{\textbf{c}}." Let us deal with a ternary finitely uniformity non-trivial theory. We apply uniformity scheme to an arbitrary \textbf{3}-predicate \textit{\textbf{C}} (I note that by applying the homogeneity scheme, we can consider any pair of elements from the triple, provided that all existential quantifiers precede universal quantifiers):
$$\forall x:C(x,x,x)\vee\neg C(x,x,x)$$
— a predicate can only be reflexive or irreflexive;
$$\forall(x,y,z):C!(x,y,z)$$
— non-triviality requires that for any triple (hereinafter — unequal!) elements there is a permutation for which the predicate is true — completeness;
$$\exists (x,y,z)\exists t:C(x,y,z)\wedge C(x,y,t)\wedge C(x,t,z)\Rightarrow\exists (x,y,z)\exists t:C(x,y,z)\wedge C(x,y,t)\wedge C(x,t,z)$$
— if there exists a t that "condenses" some triple, then the whole relation is everywhere dense; together with completeness, this means that as soon as the number of elements is greater than three (which would be one of the trivial cases), the relation is everywhere dense and the model set is infinite;
$$[\exists(x,y,z,t):C(x,y,t)\wedge C(x,t,z)\Rightarrow C(x,y,z)]\Rightarrow$$
$$\Rightarrow[\forall(x,y,z,t):C(x,y,t)]\vee C(x,t,z)\Rightarrow C(x,y,z)]$$
— an element "compacting" a pair "following" an arbitrary x will form a transitive triple with this pair; hence the whole relation is transitive;
$$[\exists(x,y,z):C(x,y,z)\Rightarrow C(z,x,y)]\Rightarrow[\forall(x,y,z):C(x,y,z)\Rightarrow C(z,x,y)]$$
— if at least one triple is cyclic, then the whole relation is cyclic; those — it is either everywhere cyclic or acyclic; in the case of cyclicity, we obtain an everywhere dense cyclic order, since the symmetry of $C(x,y,z)\wedge C(x,z,y)$ would mean the maximum completeness of the relation, i.e. — triviality of the theory; in the case of acyclicity, the \textbf{3}-order will degenerate into the usual binary ordering, which is also asymmetric (otherwise we get the most complete equivalence, which again will be trivial). We have proved the following theorem:\newline\par
\textbf{Theorem 3}: \textit{A non-trivial ternary finitely uniformity theory is equivalent to the theory of an everywhere dense total cyclic ordering, unique up to isomorphism.}\newline
(see above the same about the uniqueness of a binary order)\newline\par
The finite uniformity of an everywhere dense cyclic order (the inverse part of the theorem) is quite obvious — all axioms of such an order are uniformies constant-free assertions in which the variables are connected only by universal quantifiers. Recall also that the binary order cannot even be \textbf{1}-uniformity in the finite case. And the cyclic order is \textbf{1}-uniformity (but no more than that — \textbf{2}-uniformity is not satisfied due to the impossibility of compacting the finite set everywhere). We could also consider the \textbf{4}-order, the relation of mutual separation of the 4 points of the projective line, but the expected results are most likely similar; so let's focus on other things.\par
We have shown that the most "capacious" finitely uniformities theories are equivalent to everywhere dense linear (or cyclic) ordering, which allows such theories to describe sets of objects that are indistinguishable, but, if necessary, distinguishable. This is precisely the informational, cognitive value of such theories. We have not yet reached the number systems, for this we need algebra and the theory of homomorphisms, which will help to reveal the concepts of uniformity (and non-uniformity, which we have yet to define) from the other side. For now, I would like to draw attention to the considerable deductive power of axiom schemes, which make it possible to formulate general meanings and principles in a concentrated and productive way. However, in set theory there are already two such schemes, the importance of which can hardly be overestimated: this is the inductive scheme and the substitution scheme. All this once again indicates that mathematics needs "equality" of objects and properties of these objects, without their traditional dichotomously rigid distinction.
\section{Uniformity "at infinity"}
The question arises: how can one generalize and strengthen the concept of uniformity while remaining within the framework of relations?\par
Let's try to replace axiom schemes with axioms for "variable schemes", when not only atomic individual variables, but also sets of them are taken as such... Moreover, in the case of finite uniformity, we were actually talking about sets of elements. For convenience, we denote predicates by small Greek letters, elements by small Latin letters, and sets by capital Latin letters. The double arrow: $\Rightarrow$ denotes the derivability of formulas in our proofs. The use of sets as individual variables also leads us to the use of certain set-theoretic relations and functions. Namely, we will have Boolean operations and cardinalities of sets — they will thus constitute an atomic Boolean lattice with a measure, in the general case, a cardinal one. Unless otherwise stated, we "identify" x and {x}, not forgetting, however, that in fact: \textit{\textbf{x}} and $\{x\}$. Similar constructions are possible because $\{ \}:M\hookrightarrow2^{M}$ is an embedding of \textit{\textbf{M}} in its Boolean. So:
$$\rho(x,A)\wedge (B\subseteq A)\Rightarrow\rho(x,B)$$
$$\rho(A,x)\wedge (B\subseteq A)\Rightarrow\rho(B,x)$$
In fact, there can be any sets on the right and left of the predicate; this is equivalent to the fact that the element textit{\textbf{x}} is less than or equal to all elements of the set \textbf{\textbf{A}} (or greater than or equal, in the second sentence). We call this relation a \textit{constraint}. I am trying to introduce axioms for this relation. We go further (some version of the axiom of choice):\par\begin{center}
$Fin(A)\Rightarrow\exists! x:\rho(x,A)\wedge(x\in A)$ ~~~~ \textbf{F)}\end{center}
where \textit{\textbf{Fin(A)}} means the finiteness of \textit{\textbf{A}}. We claim:\newline\par
\textbf{Lemma 3}: \textit{From the scheme of axioms \textbf{F)} it follows that $\rho$ is a reflexive linear order on elementary individual variables \textbf{x} running through the model \textbf{M}.}\newline\par
Proof. We have:
$$\exists!x:(x\in\{x\})\Rightarrow\forall x:\rho(x,x)$$
$$\exists!x:(x\in\{y,z\})\wedge\rho(x,\{y,z\})\Longrightarrow\forall y\forall z (y\neq z):
[\rho(y,z)\vee\rho(z,y)]\wedge\neg[\rho(y,z)\wedge\rho(z,y)]$$
— the relation is reflexive, antisymmetric and linear (i.e. — complete); and now we derive transitivity:
$$\exists!t:(t\in \{x,y,z\})\wedge\rho(t,\{x,y,z\})\wedge\exists!t:(t\in \{x,y\})\wedge\rho(t,\{x,y\})\Rightarrow$$
$$\Rightarrow\forall (x,y,z):[\rho(x,y)\wedge\rho(y,z)]\rightarrow\rho(x,z)$$
Lemma proven\newline\par
So, now this order can be either a well order on a finite set or on any ordinal, or an unbounded linear order, \textbf{1}-uniformity (discrete, like integers) or finitely uniformity (everywhere dense). We also see that in any case, within any finite set, it is precisely the well order that arises, which gives rise to absolute distinguishability \textit{within} a finite set. If formula \textit{\textbf{F)}} holds for \textbf{every} set bounded below, then we have a discrete case, a \textbf{1}-uniformity but not \textbf{2}-uniformity theory in which every right segment (all elements greater than some given ) is well ordered. And if \textit{\textbf{F)}} is true \textbf{for all} sets, we get a natural ordinal.\newline\par
\textit{If, however, the theory is finitely uniformity, then a further increase in uniformity may consist in the existence of exact minimums of all right-hand segments. This turns the model set into a linearly ordered continuum.}\newline\par
This deserves more detailed consideration. A set $Q\subseteq M$ is called a \textbf{right segment} if:
$$\forall x\forall y: ((x\in Q)\wedge\rho(x,y))\Rightarrow(yin Q)$$
$$\exists t\exists x:(t\notin Q)\wedge(x\in Q)\wedge(t,x)$$
If $\rho$ is an order relation, then such a definition makes sense. And then:
$$\forall x\exists!Q_{x}:Q_{x}\equiv\{t|x\leqslant t\}$$
— will be a right segment, and the formula defining it will give an isomorphism of M into the set of all right \textit{\textbf{M}}-segments, reversely ordered by inclusion. If there is no greatest element in \textit{\textbf{M}}, then the sets $\{x<t\}$ will also be "the same" right segments. In the cases of finite linearly ordered sets, natural numbers, or discrete unbounded linear order (integers), all right segments are exhausted by this. Everywhere dense unbounded linear order (rational numbers) is more complicated. Such a theory is complete if the individual variables range only over elements of the model set. If we are talking about some of its subsets, then the theory is no longer complete. As for the right segments, some segments do not correspond to any rational numbers at all. And there are uncountable "majority" of such segments. At the same time, due to the isomorphism of the above embedding:
$$x=\bigvee\{y|y<Q_{x}\}$$
— thus, the unary predicates defining the sets $Q_{x}$ will uniquely construct an indicator for any \textit{\textbf{x}}. And the model of the theory of an everywhere dense order in the set of right segments is uniquely embedded in the model of a continually ordered set in which each right segment has a least bound (belonging to it or not). Those — any bounded set will have an exact edge. The uniformity of the continuum linear order, therefore, is higher than the uniformity of the rational straight line, because the power of the model set will become greater, and not only some, but all right segments will form a single uniformity object. And further embedding of the continuum in the set of its right segments will already be an epimorphism. In other words, all ordinally definite numbers exist. A completely similar process will also occur when considering a cyclic order on a circle (the same continuum, but only in angular or projective interpretations).\par
Consequently, \textit{\textbf{the continuum line (circle) is the most uniformity of all the sets we have considered.}} More precisely, countably homogeneous, because all this is tantamount to the existence of limits of countable ordered sequences — we are still in the "zone of action" of Archimedes' axiom. Further disclosure of the concept of uniformity will require consideration of morphisms and operations, increasing the power of models and applying methods of non-standard analysis.\par
In conclusion, I want to draw attention to the fact that objects such as $\sqrt 2$ arise quite naturally and inevitably. \textit{$\sqrt 2$ exists even in rational sets of numbers (and even in ordinary natural arithmetic), it's just that there it is not an individual, but is masked by sequences and formal sentences.} We can argue about the existence and non-existence of "infinite" objects, about the extent of their subjectivity, but they have a reality that is inseparable from the reality of finite objects. In the same way, one can start arguing about the reality of ordinary finite numbers — there are numbers, "pebbles and sticks", but what is a number? The very possibility of building a quantum computer suggests that there are no unambiguous answers even to such simple questions. The problem is not whether there are infinite constructs, but how to deal with them.
\section{NONUNIFORMITY}
First of all, we briefly outline the results of the study of uniformity.\newline\par
\textit{•	Any mathematical statement is a statement of similarity (difference) of some different (similar) objects. One possible similarity model is equivalence (its atomic form is equality; its trivial completeness is a \textbf{1}-element set). One of the possible models of difference is the relation of distinguishability, as the negation of equivalence. Distinguishability also has an independent definition: it is an irreflexive, symmetric and antitransitive (if two elements are different, then every third one will differ from at least one of them) relation:} $\forall(x,y,z):xRz\Rightarrow(yRx)\vee(yRz)$. \textit{Every process of measurement (be it length, temperature, probability, truth, etc.) is the fixation of some point on a "scale" of uniformly distinguishable points.}\newline\par
\textit{•	A mathematical theory \textbf{T} is called 1-uniformity if the following scheme is true:} $\exists x:A(x)\Rightarrow\forall x:A(x)$ \textit{for any non-constant unary formula \textbf{A(x)}. Nontrivial \textbf{1}-uniformity is in the final case possible only in the form of a well cyclic (i.e. ternary) order, to which it is equivalent. The bi-predicate has as its "minimal" \textbf{1}-uniformity realization only in both directions an infinite discrete linearly ordered set, which is order isomorphic to the set of integers. And such an order is also equivalent to the scheme of \textbf{1}-uniformity. But this order will not even be \textbf{2}-uniformity.}\newline\par
\textit{•	We call \textbf{T n}-uniformity if the theory induced by it on variables spanning all \textbf{n}-element subsets of the set spanned by \textbf{T} is \textbf{1}-uniformity. A theory that is \textbf{n}-uniformity for all finite \textbf{n} is said to be finitely uniformity. Scheme of \textbf{n}-uniformity:}
$$\exists(x_{1},...,x_{n}):\bigwedge_{1,...,n}(x_{i}\neq x_{j})\wedge A(x_{1},...,x_{n})\Longrightarrow\forall(x_{1},...,x_{n}):A!(x_{1},...,x_{n})$$
\textit{— for any \textbf{n}-ary formula \textbf{A}. A non-trivial ternary finitely uniformity theory is equivalent to the theory of everywhere dense total cyclic ordering, unique up to isomorphism. This is how the rational angles of the circle are ordered. If a theory is finitely uniformity, binary, not very simple, and nontrivial, then it is a theory of everywhere dense unbounded linear order with a single predicate. The reverse is also true. This is how the rational line is ordered. And any finite subsets of these objects will be well ordered, i.e. are completely nonuniformity.}\newline\par
\textit{•	If the theory is finitely uniformity, then a further increase in uniformity might be to have exact infimums of all right-hand segments (open intervals on the circle). This turns the model set into a linearly ordered continuum. The continuous line (circle) is the most uniformity of all the previously considered sets.}\newline\par
\textit{•	This brings us to minimality schemes. For example, F-scheme:}
$$Fin(A)\Rightarrow\exists!x:\rho(x,A)\wedge(x\in A)$$
\textit{— which is equivalent to the fact that} $\rho$ \textit{is a reflexive linear order on elementary individual variables x running through the model \textbf{M} ( \textbf{Fin(A)} means the finiteness of an arbitrary set \textbf{A} ).}\newline\par
\textbf{1}-nonuniformity meaningfully means an isomorphism of an individual domain into a set of constant-free unary formulas of the theory \textit{\textbf{T}} (moreover, each such formula is true only on the corresponding element). Let's try to formalize it. The uniqueness condition suggests that the Boolean algebra of the theory \textit{\textbf{T}} will become atomic, at least in a \textbf{1}-nonuniformity model, since the elements of a set are the atoms of its Boolean. All constant-free unary formulas \textit{\textbf{A(x)}} of the theory are Boolean algebras (we identify the equivalent ones). Let's call \textit{\textbf{A(x)}} \textbf{b}-atomic (Boolean atomic) if it is derivable in \textit{\textbf{T}}: $\exists!x:A(x)$. If this is derivable in some model \textit{\textbf{M}}, then we will call it \textbf{M}-atomic (model atomic in \textit{\textbf{M)}}. In both cases, the atomicity of \textit{\textbf{A(x)}} is equivalent to the fact that every formula with non-empty truth has a subformula with the only truth:
$$\exists x:A(x)\Rightarrow\exists!x\forall y:A_{0}(x)\wedge[A_{0}(y)\Rightarrow A(y)]$$
— in this case, $A_{0}\equiv A_{0}(x,A)$ can depend on both \textit{\textbf{A}} and \textit{\textbf{x}}. By a subformula we mean a formula from which a given \textit{\textbf{A}} follows. Such atomicity could also serve as a definition of \textbf{1}-nonuniformity. It remains only to find the explicit form of the atomic subformula.\par
$A_{0}$ must be true on a single element from the truth region of \textit{\textbf{A}}. Hence, $A_{0}$ depends on \textit{\textbf{A}} precisely through its truth region, its extension. But this is equivalent to some fixed bi-predicate $\rho(x,y)$; those —\newline\par
\textbf{Definition 4:} \textit{A theory is called \textbf{1}-nonuniformity or atomic if it contains a bi-predicate $\rho$ such that for any non-constant unary formula A there will be:}\par\begin{center}
$\exists!x\forall y:\neg A(y)\vee(A(x)\wedge\rho(x,y))$ ~~~~ \textbf{Q)}\end{center}\par
We needed some uniform procedure that matches any formula with its sub-atom. It is this role that the $\rho$-predicate performs. Thus, \textbf{Q)} is the final scheme of \textbf{M}-atomicity of theory \textit{\textbf{T}} (for the "unary" Boolean algebra of formulas). Any other schemes will be equivalent to it. Each non-empty constant-free unary formula \textit{\textbf{A(x)}} corresponds to its \textbf{M}-atomic subformula $A_{0}\equiv A_{0}(x,A)\equiv Q(x,A)$, where: \textit{\textbf{Q(x, A)}} is obtained from the scheme by deleting the existential quantifier. The fixedness of the $\rho$-predicate creates an even stronger condition: out of all possible such subformulas, the only one is chosen. It should also be noted that in the case of a \textbf{b}-atomic theory, the existence of a $\rho$-predicate naturally follows from the axioms of the theory; and in the case of \textbf{M}-atomicity, the $\rho$-predicate could also appear purely "externally", simply being added to the axioms.\par
Now it is necessary to determine the situation, in order to obtain a scheme of axioms — from the properties of the Boolean algebra \textit{\textbf{T}} to derive the properties of the elements of an individual domain — it is necessary to deal with the $\rho$-predicate. First of all, we pay attention to the fact that \textit{\textbf{Q)}} coincides in form with the uniqueness conditions for exact faces of ordered sets. We now use the condition of this uniqueness.\par
Each M-atomic subformula corresponds to the only element on which it is true. We will also prove the converse: \textbf{1}-nonuniformity of \textit{\textbf{T}} with scheme \textbf{Q)}. And since any \textbf{1}-nonuniformity \textit{\textbf{T}} must be \textbf{M}-atomic (due to the complete distinguishability of individuals), then \textbf{Q)} can be considered the definition of \textbf{1}-nonuniformity. In general, the distinguishability of individuals entails the distinguishability of any of their sequences. So, we will talk about absolute nonuniformity (as in the case of uniformity, further "building up" of nonuniformity will only be with an increase in power, either model or theoretical). In fact, we will prove even more — the equivalence of the scheme \textbf{Q)} of the theory of well order.\newline\par
\textbf{Theorem 4:} \textit{Any \textbf{1}-nonuniformity theory (i.e., a theory in which the \textbf{Q}-scheme is true) is a theory of well order. And vice versa.}\newline\par
Proof. Let there exist an element that is non-reflexive for a $\rho$-predicate; then, according to \textbf{Q)}:
$$\exists!x\forall y:\neg\rho(y,y)\Rightarrow(\neg\rho(x,x)\wedge\rho(x,y))$$\par
The resulting contradiction proves the reflexivity of the $\rho$-predicate:
$$\forall x: \rho(x,x)$$\par
Now let there be $\rho$-incomparable pairs. We define incomparability and, applying \textbf{Q)}, we get the conclusion:
$$no(x)\equiv\exists y:\neg\rho(x,y)\wedge\neg\rho(y,x)$$
$$\exists x:no(x)\Rightarrow\exists!x:[no(x)\wedge\forall y:(no(y)\Rightarrow\rho(x,y))]\Longrightarrow[no(x)\wedge\neg no(x)]$$\par
Whence we finally have the completeness of the $\rho$-predicate:
$$\forall x\forall y:\rho(x,y)\rho(y,x)$$
That is, the set of all non-reflexive elements contains its reflexive "inferior bound". Therefore, the $\rho$-predicate is reflexive:
$$\forall x:\rho(x,x)$$\par
Further, if $\exists x\exists y:\rho(x,y)\wedge\rho(x,y)\wedge(x\neq y)$; then the "minimal" of such symmetrical elements (i.e., for which there exists at least one, unequal to it, but symmetrical) cannot be unique; hence the $\rho$-predicate is antisymmetric:
$$\forall x\forall y:\rho(x,y)\wedge\rho(y,x)\Rightarrow(x=y)$$\par
Any element included in some non-transitive triple is called non-transitive. This property is expressible as a logical formula for $\rho$. In such a case, there is a single element \textit{\textbf{c}}, the "minimal" of non-transitive elements. The latter is possible in \textbf{3} cases:\par
1. ~~~~ $\exists x\exists y:\rho(c,x)\wedge\rho(x,y)\wedge\neg\rho(c,y)$ — but \textit{\textbf{y}} is also non-transitive, i.e. — $\rho(c,y)$.\par
2. ~~~~ $\exists x\exists y:\rho(y,x)\wedge\rho(x,c)\wedge\neg\rho(y,c)$ — but \textit{\textbf{x}} is also non-transitive, i.e. — $\rho(c,x)$,\par
~ ~ ~ ~ ~whence: $(x=c)$; means: $\rho(y,s)$.\par
3. ~~~~ $\exists x\exists y:\rho(y,c)\wedge\rho(c,x)\wedge\neg\rho(y,x)$\par
This has covered all cases. And they all lead to contradiction. Therefore, our predicate is transitive:
$$\forall x\forall y\forall z:\rho(x,y)\wedge\rho(y,z)\Rightarrow\rho(x,z)$$\par
We move on. If we substitute the true \textit{\textbf{A}} everywhere in the \textbf{Q}-scheme, we will see that our model set has the smallest element:
$$\exists!x\forall y:\rho(x,y)$$\par
And the last. Assume the existence of an element "above which" there is an everywhere dense segment. Here, as in all our derivations from the \textbf{Q}-scheme, it is important to find exactly the appropriate constant-free unary formula \textit{\textbf{A}}, and not to use arbitrary sets of elements in our reasoning. In fact, we have already proved that the $\rho$-predicate is a linear order with the smallest element. And now we will find out what is the density of this order. Each pair of elements can or cannot be densified by some third element lying strictly between them. And such compaction, in turn, may or may not be further compacted, "to infinity". And, if we manage to find suitable formulas for these concepts, then we can apply the \textbf{Q}-scheme to them and find the very minima that will determine the density status of the $\rho$-predicate. To shorten the derivation, we give intermediate definitions (as a result of which we arrive at the definition of the element from which the everywhere dense right segment begins):
$$P(x,y)\equiv\rho(x,y)\wedge(x\neq y)$$
$$pl(x,y)\equiv P(x,y)\wedge\exists z: P(x,z)\wedge P(z,y)$$
$$Pl(x,y)\equiv pl(x,y)\wedge\forall z:[P(x,z)\wedge P(z,y)]\Rightarrow[pl(x,z)\wedge pl(z,y)]$$
$$Pl(x)\equiv\exists y: Pl(x,y)$$\par
And now we output:
$$\exists x:Pl(x)\Rightarrow\exists!x:Pl(x)\wedge\forall y:[Pl(y)\Rightarrow\rho(x,y)]$$
$$\exists!z\exists!x:Pl(x)\wedge Pl(z)\wedge\forall y:[Pl(y)\Rightarrow\rho(x,y)]\wedge\forall t:[Pl(t)\wedge(t\neq x)\Rightarrow\rho(z,t)]$$
those — among all elements that start an everywhere dense right segment (let's call them dense on the right), there is the smallest one. But then there must also exist the smallest among all dense on the right, but not coinciding with this smallest. In this long formula (both verbally and formally) we deliberately do not resort to abbreviations (we do not write: "the element for which...", but everywhere we write: "there is only one..."), because "\textbf{Q)} — the "minimized" formulas \textit{\textbf{A}} must be constant. In fact, we have shown that the smallest of the dense ones is followed by the smallest of the dense ones that do not coincide with it. But then the "first" smallest will not be dense! We get a contradiction. The conclusion is that the order is everywhere discrete:
$$\forall x\exists!y\forall z:(x\neq y)\wedge\rho(x,y)\wedge\neg(\rho(x,z)\wedge\rho(z,y))$$
The theorem has been proven.\par
We have proved that the $rho$-predicate is an order relation. And from \textbf{Q)} follows the principle of minimality: Each set defined by a non-empty non-constant formula contains a unique smallest element. Hence our \textit{\textbf{T}} is a theory of full order. Moreover, we never used the axioms of set theory for "sets" of model elements on which the unary formulas used in the \textbf{Q}-scheme of axioms are true. Despite the fact that this would greatly simplify our proofs (sets of 2-s, 3-s of elements, with minimal ones among them, etc.). In addition, we have not used the induction scheme anywhere. We only had a \textbf{Q}-scheme with a fixed "atomizing" bi-predicate.\par
So, we have found a scheme of absolute (see above) nonuniformity; it is also a scheme of \textbf{M}-atomicity (however, it is also a certain scheme of \textbf{b}-atomicity, because the Boolean algebra of formulas of a theory with a \textbf{Q}-scheme is clearly atomic; one should not forget that this atomicity is induced by a bi-predicate, but there are also possible other ways to create an atomic \textit{\textbf{T}}), it is also a scheme of quite order. And... this is not yet a natural series, because we do not have an induction scheme. Moreover, it suffices to add the scheme of transfinite induction, and we obtain the theory of ordinals. But, and this is very important:\newline\par
\textbf{Theorem 5:} \textit{Scheme \textbf{Q)} together with the scheme of natural induction describe the natural numbers, i.e. — they are equivalent to arithmetic.}\newline\par
This is an obvious consequence of the previous reasoning.\par
Moreover, it is obvious that we can easily determine the function of "following" $(n+1)$ from these 2 schemes, for which all Peano's axioms will be fulfilled.\par
As for the scheme of induction (as well as the recursion based on it), one should not forget that for all its importance and power, it is important precisely theoretically, and precisely when we are trying to model the infinite. For practical calculations limited to some finite range of numbers, the scheme of induction always follows from other axioms and is needed as a general cognitive principle, but not as a separate axiom. The real meaning of the scheme of induction is in the statement about the cardinality of the model set (countability is the minimum of infinities). Schemes of type \textbf{Q)}, being different versions of the axiom of choice, work in all cases and always require an axiomatic formulation.\par
So, at the moment we have schemes of uniformity (in the \textbf{Definition 1)} and a \textbf{Q}-scheme of nonuniformity. Important, and not at all as trivial as it seems, is also the equality scheme:
$$\forall x:(x=x)$$
$$\forall x\forall y:(x\neq y)\vee \neg A(x)\vee A(y)$$
where \textbf{A(x)} is any formula.\par
But let's get back to nonuniformity — let's look at it from the point of view of atomicity. The \textbf{Q}-scheme, due to the strict fixation of the $\rho$-predicate, is slightly stronger than M-atomicity. Namely, textbf{M}-atomicity, in its most general form, is when the formula variable is not only \textit{\textbf{A}}, but also $\rho$ (possibly with a partial domain of definition). If we try to make $\rho$ dependent on \textit{\textbf{A}}, then any formula with full truth $A(x)\vee\neg A(x)$ will still specify a specific unique order on the entire model set \textit{\textbf{M}} — such an extension of \textit{\textbf{T}} is not essential. However, we can take $\rho$ running through an arbitrary set of complete orders on \textit{\textbf{M}}. I.e. — we will have \textit{\textbf{T}} of an arbitrary collection (or all) of completely ordered \textit{\textbf{M}}. If we have no means of distinguishing these orders, then this is the case when we know that $\exists!x:...$, but we do not know which of them elements. This is the extensional nonuniformity (or, if you like, uniformity — we have the case when the opposites almost converged), which we mentioned in the previous article. It is precisely this possibility that \textbf{M}-atomicity is wider than nonuniformity. If we distinguish all $\rho$-predicates, then again we get only an inessential extension of the theory of complete order. From this we obtain the general definition of extensional \textbf{1}-nonuniformity: $\exists!x:A(x)$ implies the undecidability of \textit{\textbf{A(x)}}, i.e. — all \textbf{M}-atomic formulas are undecidable in constants, even if constant-free closed formulas and axiom schemes are completely decidable.\par
Now let's talk about \textbf{b}-atomicity. \textbf{M}-atomicity is when any unary formula has a subformula that has the uniqueness property. This is a fruitful concept, a special case of the concept of \textbf{b}-atomicity, when any constant-free unary formula has a subformula that does not follow from any formulas that do not coincide with it. Those,— in view of the analogy with \textbf{M}-atomicity, we obtain the following definition of \textbf{b}-atomicity:\newline\par
\textbf{Definition 5:}
$$z(x)\equiv\forall y:\neg A(y)\vee[A(x)\wedge\delta(x,y)]$$\begin{center}
$\exists x:z(x)$ ~~~~ $\forall x\forall y:[z(x)\wedge z(y)]\Rightarrow[B(x)\vee\neg B(y)]$ ~~~~ \textbf{Q1)}\end{center}
\textit{— where \textbf{A, B}, are formula variables (unary and non-constant), and $\delta(x,y)$ is a fixed bi-predicate.}\newline
The second of these formulas is a scheme of \textbf{1}-uniformity in each set of "minimal" individuals (after all, we removed the uniqueness condition, but required the absence of non-empty subformulas of "minimality domains"). In general, a scheme of this type, but with an arbitrary condition \textit{\textbf{z(x)}}, can be called a scheme of conditional uniformity, and this is perhaps another way of studying gradations of uniformity.\par
Let us now consider the consequences of b-atomicity. The relation $\delta(x,y)$ on the model \textit{\textbf{M}} given by the scheme \textbf{Q1)} is a complete preordering, i.e. — \textit{\textbf{M}} is divided into well-ordered equivalence classes, whose elements are 1-uniformity within each class. This can be proved step by step, as in the proof of complete ordering for the $\rho$-predicate (then the whole proof only rephrases what has already been stated, with the replacement: "there is a unique such that..." by "there are such that..."). It is quite obvious that the intersections of all "minimality classes" of all possible non-empty formulas form the desired equivalence relation.\par
In nonuniformity \textit{\textbf{T}} this equivalence becomes an equality. And in \textbf{1}-uniformity this equivalence is trivially complete. In the general case, as we see, the \textbf{b}-atomic \textit{\textbf{T}} is a kind of transitional form. Note that in the second formula of the \textbf{Q1}-scheme we can substitute any scheme of finite, and even continual, uniformity instead of a \textbf{1}-uniformity scheme. This leads to various situations of nonuniformity. And further, in the same spirit, we can vary our schemes of axioms. In all cases, these brief schemes have great cognitive potential.\newline\par
Let's summarize what we came to in this part of the article:\newline\newline
\textit{•~~Minimality schemes immediately allow us to formulate the concept of nonuniformity. Let be any non-empty constant formula. Then the nonuniformity scheme (\textbf{Q}-scheme, which we just wrote in a slightly different form):
$$\exists!x\forall y:A(y)\Rightarrow[A(x)\wedge\rho(x,y)]$$
— each set defined by a non-empty non-constant formula contains a single "smallest" element. But the quotes can be removed, because we proved in this article that the $\rho$-predicate is a complete predicate of well-orderedness, the theory of which is equivalent to the theory of absolute (for all finite and infinite subsets of \textit{\textbf{M}}) nonuniformity. This scheme, together with the scheme of induction "from a finite" and the scheme of equality, is completely equivalent to arithmetic. And together with the scheme of transfinite induction, we obtain the theory of ordinals. The \textbf{Q}-scheme, in view of the assertion of the uniqueness of the smallest element, is also a scheme of model atomicity (it associates a non-empty formula with its minimal subformula, in the sense of implicative ordering).\newline
•~~And scheme \textbf{Q1)} treats of Boolean atomicity:
$$z(x)\equiv\forall y:\neg A(y)\vee[A(x)\wedge\delta(x,y)]$$\begin{center}
$\exists x:z(x)$ ~~~~ $\forall x\forall y:[z(x)\wedge z(y)]\Rightarrow[B(x)\vee\neg B(y)]$\end{center}
The relation $\delta(x, y)$ on the model M given by the scheme \textbf{Q1)} is a complete preordering, i.e. — \textbf{M} is divided into well-ordered equivalence classes, whose elements are \textbf{1}-uniformity within each class.}\par
This concludes the "predicate" part of the answer to the question posed at the beginning of the previous article. Namely, if we want to obtain information about an arbitrary object, we will inevitably have to find its position among other, but similar to it, continuously ordered objects that can potentially be completely ordered. And what is the role of specific numerical values and operations? The \textbf{n}-nonuniformity (nonuniformity of \textbf{n}-sequences of unequal elements) may also turn out to be not entirely trivial. This is yet to be explored.
\section{References}
\noindent [1] Joseph R. Shoenfield. Oxford.Addison-Wesley Publishing Company,\par
\textit{Mathematical Logic}, (1967).\newline
 [2] Haskell B. Curry,  McGraw-Hill Book Company, INC.\par
\textit{Foundations of Manthematical logic}, (1984).\newline
 [3] Helena Rasiowa and Roman Sikorski,  Panstwowe Wydawnlctwo Naukowe Warszawa,\par
 \textit{The Mathemayics of Metamathematics}, (1963).\newline
 [4] Alonzo Church, Princeton University Press,\par
 \textit{Introduction to Mathematical Logic}, (1956).\newline
 [5] D. Hilbert und P. Bernays. Springer-Verlag Berlin — Heidelberg — New York,\par
\textit{Grundlagen der Mathematik. I}, (1968).

\end{document}